\documentclass[12pt]{amsart}
\usepackage{amscd}      
\usepackage{amssymb}
\usepackage{amsmath}
\usepackage{xypic}      
\LaTeXdiagrams          
\usepackage[all,v2]{xy}
\xyoption{2cell} \UseAllTwocells \xyoption{frame} \CompileMatrices
\allowdisplaybreaks[3]

\usepackage{latexsym}
\usepackage{epsfig}
\usepackage{amsfonts}
\usepackage{enumerate}
\usepackage{times}

\newtheorem{prop}{Proposition}[section]

\newtheorem{thm}[prop]{Theorem}

            {\nolinebreak $\Box$ \end{trivlist}}

\newcommand{\noprint}[1]{}

\newcommand{\zz}{{\mathbb Z}}

\newcommand{\cc}{{\mathbb C}}

\newcommand{\ldiag}[1]%
       {\makebox[0cm]{${\scriptstyle#1}\downarrow\phantom{\scriptstyle#1}$}}
\newcommand{\ldiagup}[1]%
       {\makebox[0cm]{${\scriptstyle#1}\uparrow\phantom{\scriptstyle#1}$}}
\newcommand{\rdiag}[1]%
       {\makebox[0cm]{$\phantom{\scriptstyle#1}\downarrow{\scriptstyle#1}$}}
\newcommand{\sediagr}[1]%
       {\makebox[0cm]{$\phantom{\scriptstyle#1}\searrow{\scriptstyle#1}$}}
\newcommand{\nediagr}[1]%
       {\makebox[0cm]{$\phantom{\scriptstyle#1}\nearrow{\scriptstyle#1}$}}
\newcommand{\rdiagup}[1]%
       {\makebox[0cm]{$\phantom{\scriptstyle#1}\uparrow{\scriptstyle#1}$}}
\newcommand{\swdiag}[1]%
       {\makebox[0cm]{$\phantom{\scriptstyle#1}\swarrow{\scriptstyle#1}$}}
\newcommand{\sediag}[1]%
       {\makebox[0cm]{${\scriptstyle#1}\searrow\phantom{\scriptstyle#1}$}}
\newcommand{\nediag}[1]%
       {\makebox[0cm]{${\scriptstyle#1}\nearrow\phantom{\scriptstyle#1}$}}

\newcommand{\doublearrowstack}[2]%
                      {{{{\scriptstyle#1}\atop{\textstyle\longrightarrow}}\atop{{\textstyle\longrightarrow}\atop{\scriptstyle#2}}}}
\newcommand{\rightleftarrowstack}[2]%
                      {{{{\scriptstyle#1}\atop{\textstyle\longrightarrow}}\atop{{\textstyle\longleftarrow}\atop{\scriptstyle#2}}}}
\newcommand{\leftrightarrowstack}[2]%
                      {{{{\scriptstyle#1}\atop{\textstyle\longleftarrow}}\atop{{\textstyle\longrightarrow}\atop{\scriptstyle#2}}}}

\newcommand{\overtoparrow}%
{\makebox[0cm]{\beginpicture \setcoordinatesystem units
<.8cm,.4cm> point at 0 0 \setplotarea x from -3 to 3, y from 0 to
1 \setquadratic \plot -3 0 0 1 3 0 / \put{\vector(3,-1){0}}[Bl] at
3 0
\endpicture}}

\newcommand{\underbottomarrow}%
{\makebox[0cm]{\beginpicture \setcoordinatesystem units
<.8cm,.4cm> point at 0 0 \setplotarea x from -3 to 3, y from 0 to
1 \setquadratic \plot -3 1 0 0 3 1 / \put{\vector(3,1){0}}[Bl] at
3 1
\endpicture}}

\newcommand{\ses}[5]%
{0\longrightarrow#1\stackrel{#2}{ \longrightarrow}#3\stackrel{#4}{
\longrightarrow}#5\longrightarrow0}

\newcommand{\dt}[6]%
{#1\stackrel{#2}{longrightarrow}#3
\stackrel{#4}{\longrightarrow}#5 \stackrel{#6}{\longrightarrow}
#1[1]}

\newcommand{\cat}[1]%
{(\mbox{\rm #1})}

\def\Label#1{\label{#1}{\tt [#1]}\phantom{h}}
\def\Label{\label}

\theoremstyle{definition}

\theoremstyle{remark}

\theoremstyle{remark}

\numberwithin{equation}{section}

\newcommand{\Mbar}{\overline{\M}}

\newcommand{\X}{\mathcal{X}}

\newcommand{\M}{\mathcal{M}}
\newcommand{\C}{\mathcal{C}}

\newcommand{\B}{\mathcal{B}}

\def\<{\left\langle}
\def\>{\right\rangle}

\title{Virasoro constraints and descendant Hurwitz-Hodge Integrals}
\author{Yunfeng Jiang}
\address{Department of Mathematics\\ University of Utah\\ 155 South 1400 East JWB233 \\Salt Lake City\\ UT 84112\\ USA}
\email{jiangyf@math.utah.edu}
\author{Hsian-Hua Tseng}
\address{Department of Mathematics\\ University of Wisconsin-Madison\\ Van Vleck Hall, 480 Lincoln Drive \\Madison\\ WI 53706-1388 \\ USA}
\email{tseng@math.wisc.edu}
\date{\today}

\begin{document}
\begin{abstract}
Virasoro constraints are applied to degree zero Gromov-Witten theory of weighted projective stacks $\mathbb{P}(1,N)$ and $\mathbb{P}(1,1,N)$ to obtain formulas of descendant cyclic Hurwitz-Hodge integrals in higher genera.
\end{abstract}

\maketitle

\section{Introduction}

In \cite{JT} the general theory of Frobenius manifolds and Virasoro constraints (see e.g. \cite{DZ}, \cite{gi_gwi}, \cite{gi_frob}) is applied to write down Virasoro constraints for orbifold Gromov-Witten theory. Following the idea\footnote{This was suggested to the second author by E. Getzler in May 2004.} in the work \cite{GP}, Virasoro constraints are also applied to degree zero Gromov-Witten theory to deduce formulas for genus zero Hurwitz-Hodge integrals. In this paper we carry out the same idea in higher genus, that is, we derive formulas for higher genus descendant cyclic Hurwitz-Hodge integrals by applying Virasoro constraints to weighted projective stacks. Our main results, Theorems \ref{curvevir-1}--\ref{curvevir} and Theorems \ref{surfacevir-1}--\ref{surfacevir}, may be viewed as generalizations of the famous $\lambda_g$-theorem and Faber's intersection number conjecture. 

We closely follow the notations used in \cite{JT}. The readers are referred to \cite{JT} and references therein for notations, conventions, and general discussions on Gromov-Witten theory for Deligne-Mumford stacks and their Virasoro constraints. Stable maps to weighted projective stacks are discussed briefly in Section \ref{degree_0_maps}. The applications of Virasoro constraints to degree zero theory of weighted projective lines and planes are carried out in Sections \ref{weighted_P1} and \ref{weighted_P2} respectively. A brief discussion on the dimension three case is given in Section \ref{threefolds}.

\subsection*{Acknowledgments}
Both authors thank the organizers of the workshop {\em Recent progress of the moduli of curves} in Banff International Research Station, where part of this work was pursued. 

\section{Degree zero twisted stable maps to weighted projective spaces}\label{degree_0_maps}
In this section basics about degree zero twisted stable maps to weighted projective spaces will be discussed. 

Fix an integer $N>1$. Let $$\mathcal{X}=\mathbb{P}(\underbrace{1,\cdots,1}_{d},N)$$ be the $d$-dimensional weighted projective stack with weights $(1,1,\cdots,1,N)$. There is a unique stacky point $p=[0,\cdots,0,1]$ on $\X$ with isotropy group the cyclic group $\mathbb{Z}_{N}$. We shall identify $\mathbb{Z}_N$ with the group of $N$-th roots of unity. Put $\omega:=\exp(\frac{2\pi\sqrt{-1}}{N})$.

Fix $g\geq 1$ throughout. Let $f: \mathcal{C}\to \mathcal{X}$ be a degree zero genus $g$ twisted stable map and 
\begin{equation}\label{monodromy_rep}
\pi_1^{orb}(\C)\to \pi_1^{orb}(\X),
\end{equation}
be the associated monodromy representation. If (\ref{monodromy_rep}) is nontrivial, then $f$ must factor through $p\simeq \mathcal{B}\zz_N$. By representability, this happens when there are stacky points on $\C$. Suppose the stack structure on the domain curve $\mathcal{C}$ is given by the following tuple of elements in $\zz_N$ (here $n_{1},\cdots,n_{N-1}\in \mathbb{Z}_{> 0}$),
\begin{equation}\label{ramification}
\mathbf{x}=(\underbrace{\omega,\cdots,\omega}_{n_{1}},\underbrace{\omega^{2},\cdots,\omega^{2}}_{n_{2}},\cdots,\underbrace{\omega^{i},\cdots,\omega^{i}}_{n_{i}},\cdots,\underbrace{\omega^{N-1},\cdots,\omega^{N-1}}_{n_{N-1}}).
\end{equation}

Let $C$ be the coarse curve of $\C$. The morphism $f$ is equivalent to an admissible $\mathbb{Z}_{N}$-cover $\widetilde{\mathcal{C}}\to C$ with monodromies specified by $\mathbf{x}$. Let $h$ denote the genus of $\widetilde{\C}$.

Let $\overline{\mathcal{M}}_{g,n+\Sigma_{i}n_{i}}(\mathcal{X},0,\mathbf{x})$ be
the moduli stack of degree zero genus $g$ twisted stable maps $f:\C\to \X$ to $\X$ with
$n$ non-stacky marked points and $\Sigma_{i}n_{i}$ stacky marked points of type $\mathbf{x}$. Let 
$$q: \overline{\mathcal{M}}_{n+\Sigma_{i}n_{i}}(\mathcal{X},0,\mathbf{x})\to \overline{\M}_{g,n+\Sigma_{i}n_{i}}$$
be the forgetful map to the moduli space of genus $g$ stable $(n+\Sigma_{i}n_{i})$-marked curves. Descendant classes on $\overline{\mathcal{M}}_{n+\Sigma_{i}n_{i}}(\mathcal{X},0,\mathbf{x})$
are defined via pull-back, $$\psi_i:=q^*\psi_i.$$  Degree zero descendent orbifold Gromov-Witten invariants of $\X$ are defined by
\begin{equation}\Label{orbinv}
\langle\tau_{l_{1}}\cdots\tau_{l_{n}}\widetilde{\tau}_{k_{1}}\cdots\widetilde{\tau}_{k_{\Sigma_{i}
n_{i}}}\rangle_{g}^{\mathcal{X}}
:=\int_{[\overline{\mathcal{M}}_{g,n+\Sigma_{i}n_{i}}(\mathcal{X},0,\mathbf{x})]^{vir}}
\psi^{l_{1}}\cdots\psi^{l_{n}}\psi^{k_{1}}\cdots\psi^{k_{\Sigma_{i}
n_{i}}}.
\end{equation}
Here $\tau$ indicates descendant insertion from non-stacky marked points, and $\widetilde{\tau}$ indicates descendant insertion from stacky marked points.

Let $L_\omega$ denote the line bundle over $\mathcal{B}\zz_N$ defined by the $\zz_N$-representation $\cc$ on which $\omega\in \zz_N$ acts by multiplication by $\omega$. It can be seen from the definition of $\X$ that the normal bundle of $\mathcal{B}\mathbb{Z}_{N}\subset \mathcal{X}$ is $L_{\omega}^{\oplus d}$. Let $\Mbar_{g,n+\sum_i n_i}(\B\zz_N, \mathbf{x})$ denote the moduli stack of genus $g$ stable maps to $\B\zz_N$ with marked points of stack type $\mathbf{x}$. Consider the universal diagram associated to stable maps to $\B\zz_N$ of stack type $\mathbf{x}$:
$$\xymatrix{
~\widetilde{\mathcal{C}}\dto_{p}\rto & pt\dto^{} \\
~\mathcal{C}\rto^{f}\dto_{\pi} &  ~\mathcal{B}\mathbb{Z}_{N}\\
\overline{\mathcal{M}}_{g,n+\Sigma_{i}
n_{i}}(\mathcal{B}\mathbb{Z}_{N}, \mathbf{x}).&~}$$ 

Let $\widetilde{\pi}=\pi\circ p: \widetilde{\mathcal{C}}\to \overline{\mathcal{M}}_{n+\Sigma_{i}n_{i}}(\mathcal{B}\mathbb{Z}_{N},\mathbf{x})$ be the composite map. The $\zz_N$-action on $\widetilde{\C}$ induces an action on the dual of the Hodge bundle $\mathbb{E}^{\vee}=R^{1}\widetilde{\pi}_{*}\mathcal{O}$. This yields the decomposition of $\mathbb{E}^\vee$ into $\omega$-eigenbundles:
$$\mathbb{E}^{\vee}=\mathbb{E}_{1}^{\vee}\oplus \mathbb{E}_{\omega}^{\vee}
\oplus\cdots\oplus \mathbb{E}_{\omega^{N-1}}^{\vee}.$$
Following the convention in \cite{BGP}, $\mathbb{E}_{\omega^i}^\vee$ is the eigen-bundle on which 
$\omega$ acts with eigenvalue $\omega^i$.

It is easy to see that
$$R^{1}\pi_{*}f^{*}(L_{\omega})=\mathbb{E}_{\omega^{N-1}}^{\vee}\simeq \mathbb{E}_{\omega}.$$
An analysis on obstruction theory yields
\begin{align}\Label{orbinv2}
\langle\tau_{l_{1}}\cdots\tau_{l_{n}}\widetilde{\tau}_{k_{1}}\cdots\widetilde{\tau}_{k_{\Sigma_{i}
n_{i}}}\rangle_{g}^{\mathcal{X}}
&=\int_{\overline{\mathcal{M}}_{g,n+\Sigma_{i}n_{i}}(\mathcal{B}\mathbb{Z}_{N}, \mathbf{x})}
\psi^{l_{1}}\cdots\psi^{l_{n}}\psi^{k_{1}}\cdots\psi^{k_{\Sigma_{i}
n_{i}}} e(\mathbb{E}_{\omega^{N-1}}^{\vee}\oplus\cdots\oplus
\mathbb{E}_{\omega^{N-1}}^{\vee})\nonumber \\
&=\int_{\overline{\mathcal{M}}_{g,n+\Sigma_{i}n_{i}}(\mathcal{B}\mathbb{Z}_{N}, \mathbf{x})}
\psi^{l_{1}}\cdots\psi^{l_{n}}\psi^{k_{1}}\cdots\psi^{k_{\Sigma_{i}
n_{i}}} \lambda_{r_{1}}^{d},
\end{align}
where $e$ is the  Euler class and $\lambda_{r_{1}}$ the top Chern class of $\mathbb{E}_{\omega}$. The ranks $r_{1}=rank(\mathbb{E}_{\omega})$, $r_{N-1}=rank(\mathbb{E}_{\omega^{N-1}})$ are easily calculated by Riemann-Roch, 
$$r_{1}=\sum_{i=1}^{N-1}n_i\frac{i}{N}+g -1, \quad r_{N-1}=\sum_{i=1}^{N-1}n_i\frac{N-i}{N}+g-1.$$
Also note that $r_{1}+r_{N-1}-1=2g+\Sigma_{i}n_{i}-3$.

In what follows we consider Hurwitz-Hodge integrals (\ref{orbinv2}) which arise from $\mathbb{P}(1,N)$ and $\mathbb{P}(1,1,N)$. We introduce some notations for such integrals.

\begin{equation}\Label{orbinv3}
\langle\tau_{l_{1}}\cdots\tau_{l_{n}}\widetilde{\tau}_{k_{1}}\cdots\widetilde{\tau}_{k_{\Sigma_{i}
n_{i}}}|\lambda_{r_{1}}\rangle^{h}_{g}
=\int_{\overline{\mathcal{M}}_{g,n+\Sigma_{i}n_{i}}(\mathcal{B}\mathbb{Z}_{N},\mathbf{x})}
\psi^{l_{1}}\cdots\psi^{l_{n}}\psi^{k_{1}}\cdots\psi^{k_{\Sigma_{i}
n_{i}}} \lambda_{r_{1}},
\end{equation}

\begin{equation}\Label{orbinv4}
\langle\tau_{l_{1}}\cdots\tau_{l_{n}}\widetilde{\tau}_{k_{1}}\cdots\widetilde{\tau}_{k_{\Sigma_{i}
n_{i}}}|\lambda^{2}_{r_{1}}\rangle^{h}_{g}
=\int_{\overline{\mathcal{M}}_{g,n+\Sigma_{i}n_{i}}(\mathcal{B}\mathbb{Z}_{N}, \mathbf{x})}
\psi^{l_{1}}\cdots\psi^{l_{n}}\psi^{k_{1}}\cdots\psi^{k_{\Sigma_{i}
n_{i}}} \lambda^{2}_{r_{1}}
\end{equation}
and
\begin{multline}\label{orbinv5}
\langle\langle\tau_{l_{1}}\cdots\tau_{l_{n}}\widetilde{\tau}_{k_{1}}\cdots\widetilde{\tau}_{k_{\Sigma_{i}
n_{i}}}|\lambda_{r_{1}}\rangle\rangle^{h}_{g}=\\ \sum_{M\geq
0}\frac{1}{M!}\sum_{b_1,\cdots,b_M\geq 0}t_{b_{1}}\cdots
t_{b_{M}}\langle\tau_{b_{1}}\cdots\tau_{b_{M}}\tau_{l_{1}}\cdots\tau_{l_{n}}\widetilde{\tau}_{k_{1}}\cdots\widetilde{\tau}_{k_{\Sigma_{i}
n_{i}}}|\lambda_{r_{1}}\rangle^{h}_{g},
\end{multline}

\begin{multline}\label{orbinv6}
\langle\langle\tau_{l_{1}}\cdots\tau_{l_{n}}\widetilde{\tau}_{k_{1}}\cdots\widetilde{\tau}_{k_{\Sigma_{i}
n_{i}}}|\lambda^{2}_{r_{1}}\rangle\rangle^{h}_{g}=\\ \sum_{M\geq
0}\frac{1}{M!}\sum_{b_1,\cdots,b_M\geq 0}t_{b_{1}}\cdots
t_{b_{M}}\langle\tau_{b_{1}}\cdots\tau_{b_{M}}\tau_{l_{1}}\cdots\tau_{l_{n}}\widetilde{\tau}_{k_{1}}\cdots\widetilde{\tau}_{k_{\Sigma_{i}
n_{i}}}|\lambda^{2}_{r_{1}}\rangle^{h}_{g}.
\end{multline}

The following convention will be used. The symbols 
$$\<\star|\lambda_{r_1}\>_g^h, \<\star|\lambda_{r_1}^2\>_g^h$$ 
where only {\em non-stacky} insertions occur in $\star$ denote integrals over components of the relevant moduli spaces which parametrize maps with non-trivial monodromy representations. In the absence of cover genus $h$, the symbols 
$$\<\star|\lambda_g\>_g, \<\star|\lambda_g\lambda_{g-1}\>_g$$
denotes integrals over the usual moduli spaces of pointed curves. A similar convention is imposed for double bracket notations.

\section{Virasoro constraints for weighted projective lines}\label{weighted_P1}

Let $\mathcal{X}=\mathbb{P}(1,N)$. The Chen-Ruan orbifold
cohomology $H^{*}_{CR}(\mathbb{P}(1,N))$ has the following natural 
generators: $$1\in H^{0}_{CR}(\mathbb{P}(1,N)), \,
\xi \in H^{2}_{CR}(\mathbb{P}(1,N)), \,
\gamma_{j}\in H^{\frac{2j}{N}}_{CR}(\mathbb{P}(1,N)), \,
\text{ for }1\leq j\leq N-1.$$

We will use the following coordinates for the corresponding descendants:
\begin{center}
\begin{tabular}{lp{0.75\textwidth}}
coordinate & \quad descendant\\
$t_i, \, i\geq 0$ & \quad $\tau_i(1)$\\
$s_{i}, \, i\geq 0$ & \quad $\tau_{i}(\xi)$\\
$\alpha^{j}_{i}, \, i\geq 0$ & \quad $\tau_{i}(\gamma_{j})$ for $1\leq j\leq N-1$.
\end{tabular}
\end{center}

Let $\mathcal{F}_{\mathcal{X}}^{g,0}$ be the degree zero genus $g$
orbifold Gromov-Witten potential of $\mathcal{X}$. Let 
$$\mathcal{D}_{\mathcal{X}}^{0}=\exp\left(\sum_{g\geq 0}\hbar^{g-1}\mathcal{F}_{\mathcal{X}}^{g,0}\right).$$ 

For $1\leq i\leq N-1$ let 
\begin{equation}\label{notation}
M_{i}:=\sum_{a=0}^{i}n_{a},
\end{equation} 
where $n_{0}:=0$. This notation is used in Section \ref{weighted_P1} and \ref{weighted_P2}.

\begin{thm}
We have
\begin{align*}\frac{\mathcal{L}_{k}\mathcal{D}_{\mathcal{X}}^{0}}{\mathcal{D}_{\mathcal{X}}^{0}}=
\sum_{g=0}^{\infty}\hbar^{g-1}(-1)^{g}\Bigl((-2)x_{g}^{k}(\mathbf{t})+\sum_{l=0}^{\infty}s_{l}\cdot
y_{g,l}^{k}(\mathbf{t}) +\sum_{h=0}^{\infty}
w^{h,k}_{g,l}(\mathbf{t}) \\
+\sum_{h=0}^{\infty}~\sum_{k_{1},\cdots,k_{\Sigma
n_{i}}> 0}\alpha_{k_{1}}^{1}\cdots\alpha_{k_{\Sigma_{i}n_{i}}}^{N-1}\cdot
z_{h;k_{1},\cdots,k_{\Sigma_{i}n_{i}}}^{k}(\mathbf{t})\Bigr),
\end{align*} where

\begin{multline*}
x_{g}^{k}(\mathbf{t})=
-[1]_{0}^{k}\langle\langle\tau_{k+1}|\lambda_{g-1}\rangle\rangle_{g}
+\sum_{m=1}^{\infty}[m]_{0}^{k}t_{m}\langle\langle\tau_{k+m}|\lambda_{g-1}\rangle\rangle_{g}
+ [1]_{1}^{k}\langle\langle\tau_{k}|\lambda_{g}\rangle\rangle_{g} \\
-\sum_{m=0}^{\infty}[m]_{1}^{k}t_{m}\langle\langle\tau_{k+m-1}|\lambda_{g}\rangle\rangle_{g}
-\frac{1}{2}\sum_{m=0}^{k-2}\sum_{g=g_{1}+g_{2}}(-1)^{m+1}[-m-1]_{1}^{k}\langle\langle\tau_{m}|\lambda_{g_{1}}\rangle\rangle_{g_{1}}
\langle\langle\tau_{k-m-2}|\lambda_{g_{2}}\rangle\rangle_{g_{2}},
\end{multline*}
\begin{equation*}
y_{g,l}^{k}(\mathbf{t})=
-[1]_{0}^{k}\langle\langle\tau_{k+1}\tau_{l}|\lambda_{g}\rangle\rangle_{g}
+\sum_{m=1}^{\infty}[m]_{0}^{k}t_{m}\langle\langle\tau_{k+m}\tau_{l}|\lambda_{g}\rangle\rangle_{g}
+
[l+1]_{0}^{k}\langle\langle\tau_{k+l}|\lambda_{g}\rangle\rangle_{g}
\end{equation*}

$$w_{g,l}^{h,k}(\mathbf{t})=-[1]_{0}^{k}\langle\langle\tau_{k+1}|\lambda_{r_{1}}\rangle\rangle_{g}^{h}+\sum_{m=1}^{\infty}[m]_{0}^{k}t_{m}
\langle\langle\tau_{k+m}|\lambda_{r_{1}}\rangle\rangle_{g}^{h}$$
and
\begin{multline*}
z_{h;k_{1},\cdots,k_{\Sigma_{i}n_{i}}}^{k}(\mathbf{t})=
-[1]_{0}^{k}\langle\langle\tau_{k+1}\widetilde{\tau}_{k_{1}}\cdots
\widetilde{\tau}_{k_{\Sigma_{i}n_{i}}}|\lambda_{r_{1}}\rangle\rangle^{h}_{g}
+\sum_{m=1}^{\infty}[m]_{0}^{k}t_{m}\langle\langle\tau_{k+m}\widetilde{\tau}_{k_{1}}\cdots \widetilde{\tau}_{k_{\Sigma_{i}n_{i}}}|\lambda_{r_{1}}\rangle\rangle^{h}_{g}\\
+ \sum_{i=1}^{N-1}\sum_{j=M_{i-1}+1}^{M_{i}}[k_{j}+\tfrac{i}{N}]_{0}^{k}
\langle\langle\widetilde{\tau}_{k_{1}}\cdots\widetilde{\tau}_{k_{n_{i-1}}}\widetilde{\tau}_{k_{n_{i-1}+1}}\cdots\widetilde{\tau}_{k+k_{j}}\cdots\widetilde{\tau}_{k_{n_{i}}}\widetilde{\tau}_{k_{n_{i}+1}}
\cdots\widetilde{\tau}_{k_{\Sigma_{i}n_{i}}}|\lambda_{r_{1}}\rangle\rangle^{h}_{g}.
\end{multline*}
\end{thm}

\begin{proof}
The Virasoro operator $\mathcal{L}_{k}$, $k>0$ is given by
\begin{multline*}
\mathcal{L}_{k}= -[1]_{0}^{k}\partial_{t_{k+1}}+\sum_{m=0}^{\infty}\Bigl([m]_{0}^{k}t_{m}\partial_{t_{k+m}}
+[m+1]_{0}^{k}s_{m}\partial_{s_{k+m}}+\sum_{i=1}^{N-1}[m+\tfrac{i}{N}]_{0}^{k}\alpha_{m}^{i}\partial_{\alpha_{k+m}^{i}}\Bigr)\\ +
2\Bigl(-[1]_{1}^{k}\partial_{s_{k}}+\sum_{m=0}^{\infty}[m]_{1}^{k}t_{m}\partial_{s_{k+m-1}} 
+\frac{\hbar}{2}\sum_{m=0}^{k-2}(-1)^{m+1}[-m-1]_{1}^{k}\partial_{s_{m}}\partial_{s_{k-m-2}}\Bigr).
\end{multline*}
As discussed in the previous section, contributions to Gromov-Witten invariants from components of stable map spaces which parametrize maps with non-trivial monodromies (\ref{monodromy_rep}) are given as Hurwitz-Hodge integrals (\ref{orbinv2}). Components which parametrize maps with trivial monodromy representation only occurs in case when marked points are all non-stacky. A component that parametrizes degree $0$ genus $g$ $n$-pointed stable maps to $\mathbb{P}(1,N)$ with trivial monodromy is isomorphic to a product $\overline{\M}_{g,n}\times \mathbb{P}(1,N)$.
From \cite{GP}, the obstruction bundle in this case is $T_{\mathbb{P}(1,N)}\boxtimes \mathbb{E}^{\vee}$, where $\mathbb{E}$ is the usual Hodge bundle. Contributions from such a component are thus integrals against the Euler class
$$(-1)^{g}e(T_{\mathbb{P}(1,N)}\boxtimes \mathbb{E}^{\vee})=\lambda_{g}-c_{1}(\mathbb{P}(1,N))\lambda_{g-1}.$$

So the degree zero orbifold Gromov-Witten potential of $\mathbb{P}(1,N)$ is
\begin{multline*}
\mathcal{D}_{\mathcal{X}}^{0}=
\exp\Bigl((-2)\sum_{g=1}^{\infty}(-1)^{g}\hbar^{g-1}\langle\langle |\lambda_{g-1}\rangle\rangle_{g}
+\sum_{g=0}^{\infty}(-1)^{g}\hbar^{g-1}\Bigl(\sum_{m=0}^{\infty}s_{m}\langle\langle\tau_{m}|\lambda_{g}\rangle\rangle_{g}\Bigr) \\
+\sum_{g=0}^{\infty}(-1)^{g}\hbar^{g-1}\sum_{h=0}^{\infty}
\langle\langle~ |\lambda_{r_{1}}\rangle\rangle^{h}_{g}
+\sum_{g=0}^{\infty}(-1)^{g}\hbar^{g-1}\sum_{h=0}^{\infty}\sum_{k_{1},\cdots,k_{\Sigma_{i}n_{i}}> 0}\alpha^{1}_{k_{1}}\cdots
\alpha^{N-1}_{k_{\Sigma_{i}n_{i}}}\langle\langle\widetilde{\tau}_{k_{1}}\cdots
\widetilde{\tau}_{k_{\Sigma_{i}n_{i}}}|\lambda_{r_{1}}\rangle\rangle^{h}_{g}\Bigr).
\end{multline*}
Applying the operator $\mathcal{L}_{k}$ to
$\mathcal{D}_{\mathcal{X}}^{0}$ we obtain the result.
\end{proof}

It follows that degree zero Virasoro constraints for $\mathbb{P}(1,N)$ are equivalent to the vanishing of $x_{g}^{k}(\mathbf{t}), y_{g,l}^{k}(\mathbf{t})$,
$w_{g,l}^{h,k}(\mathbf{t})$ and $z_{h,k_{1},\cdots,k_{\Sigma_{i}n_{i}}}^{k}(\mathbf{t})$.

\subsubsection{\textbf{Vanishing of $w_{g,l}^{h,k}(\mathbf{t})$}.}

\begin{thm}\Label{curvevir-1}
We have
$$
\langle\tau_{l_{1}}
\cdots\tau_{l_{n}}|\lambda_{r_{1}}\rangle^{h}_{g}=
\binom{2g+n-2}{l_1,\ldots, l_n}\langle\tau_{2g-1}|\lambda_{r_{1}}\rangle_g^h. 
$$
\end{thm}

\begin{proof}
The vanishing of
$w_{g,l}^{h,k}(\mathbf{t})$ gives the vanishing of its Taylor coefficients, $$\tfrac{1}{[1]_{0}^{k}}\partial_{t_{l_{1}}}\cdots \partial_{t_{l_{n}}}w_{g,l}^{h,k}(0)=0.$$
An explicit calculation shows that $\tfrac{1}{[1]_{0}^{k}}\partial_{t_{l_{1}}}\cdots \partial_{t_{l_{n}}}w_{g,l}^{h,k}(0)$ is the right side of the following
\begin{equation}\Label{recursion-1}
0=-\langle\tau_{k+1}\tau_{l_{1}}
\cdots\tau_{l_{n}}|\lambda_{r_{1}}\rangle^{h}_{g}+\sum_{i=1}^{n}\tfrac{(l_{i}+k)!}{(l_{i}-1)!(k+1)!}\langle\tau_{l_{1}}
\cdots\tau_{l_{i}+k}\cdots\tau_{l_{n}}|\lambda_{r_{1}}\rangle^{h}_{g}
\end{equation}
The result follows by solving the recursion (\ref{recursion-1}).
\end{proof}

Let $\Mbar_{g,1}(\B\zz_N)$ be the moduli stack of genus $g$ stable maps to $\B\zz_N$ with one non-stacky marked point. The initial values in Theorem \ref{curvevir-1} are integrals over components of $\Mbar_{g,1}(\B\zz_N)$ which parametrize maps with non-trivial monodromy representations. Those integrals may be computed from Hurwitz-Hodge integrals $\int_{\overline{\mathcal{M}}_{g,1}(\mathcal{B}\mathbb{Z}_{N})}\lambda_{r_{1}}\psi^{2g-1}$ by subtracting contributions from components which parametrize maps with trivial monodromy. Such contributions are easily found to be 
$$\frac{1}{N}\int_{\Mbar_{g,1}}\lambda_{g-1}\psi^{2g-1},$$
which has been computed in \cite{FP1}. The answer may be given in generating series as 
\begin{equation}\label{1-pt-hodge}
1+\sum_{g>0}\sum_{l=0}^g t^{2g}z^l\int_{\Mbar_{g,1}}\psi_1^{2g-2+l}\lambda_{g-l}=\left(\frac{t/2}{\sin(t/2)}\right)^{z+1}.
\end{equation} 

Hurwitz-Hodge integrals of the form
$$\int_{\Mbar_{g,1}(\mathcal{B}\mathbb{Z}_{N})}\lambda_{r_{1}-i}\psi^{2g-1+i}$$
for $i\geq 0$ are computed in \cite{JPT} by means of an ELSV-type formula and exact evaluation of double Hurwitz numbers \cite{GJV}. The answer may be given in generating series as 
\begin{equation}\label{1-pt-h-hodge}
\frac{1}{N}+\sum_{g>0}\sum_{l=0}^g t^{2g}z^l\int_{\Mbar_{g,1}(\B\zz_N)}\psi_1^{2g-2+l}\lambda_{g-l}=\frac{1}{N}\left(\frac{Nt/2}{\sin(Nt/2)}\right)^z\frac{t/2}{\sin(t/2)}.
\end{equation}
Therefore (note that $r_1=g-1$)
\begin{equation}\label{1-pt-initial}
\sum_{g>0}\sum_{l=0}^g t^{2g}z^l\<\tau_{2g-2+l}|\lambda_{r_1+1-l}\>_g^h=\frac{1}{N}\frac{t/2}{\sin(t/2)}\left(\left(\frac{Nt/2}{\sin(Nt/2)}\right)^z-\left(\frac{t/2}{\sin(t/2)}\right)^z\right),
\end{equation}
and the initial values occur in the coefficient of the $z$ term.

\subsubsection{\textbf{Vanishing of $z_{h,k_{1},\cdots,k_{\Sigma_{i}n_{i}}}^{k}(\mathbf{t})$}.}

Now let
\begin{equation}\Label{constant1}
\Gamma_{j,g}=\langle\widetilde{\tau}_{a}|\lambda_{r_{1}}\rangle^{h}_{g}=\int_{\Mbar_{g,\sum n_{i}}(\mathcal{B}\mathbb{Z}_{N})}\psi^{a}\lambda_{r_{1}}, \quad \text{where } 
a:=2g-2+\sum_{i=1}^{N-1}n_{i}-\sum_{i=1}^{N-1}\frac{i}{N}n_{i}.
\end{equation}

Write $\mathbf{\Gamma_{g}}:=(\Gamma_{j,g})_{1\leq j\leq \Sigma_{i}n_{i}}$ as a column vector. Let $\mathbf{c_{g}}:=(c_{j,g})_{1\leq j\leq \Sigma_{i}n_{i}}$ be another column vector. Let the index $i$ vary from $1$ to $N-1$ and define a $\Sigma_{i}n_{i}\times\Sigma_{i}n_{i}$ square matrix $A=(a_{st})$ by

\begin{equation}\Label{matrix}
a_{st}:=\begin{cases}\frac{i}{N}+a&\text{if
}~ M_{i-1}<s=t\leq M_{i}\,;\\
\frac{i}{N}&\text{if}~M_{i-1}<t\leq M_{i}~\text{and}~s\neq t\,.\end{cases}
\end{equation}

Alternatively,
$$
A=\left[
\begin{array}{cccccccccc}
  \tfrac{1}{N}+a&\cdots&\tfrac{1}{N}&\tfrac{2}{N}&\cdots&\tfrac{2}{N}&\cdots&\tfrac{N-1}{N}&\cdots&\tfrac{N-1}{N} \\
  \tfrac{1}{N}&\cdots&\tfrac{1}{N}&\tfrac{2}{N}&\cdots&\tfrac{2}{N}&\cdots&\tfrac{N-1}{N}&\cdots&\tfrac{N-1}{N} \\
  \cdots&\cdots&\cdots&\cdots&\cdots&\cdots&\cdots&\cdots&\cdots&\cdots\\
  \tfrac{1}{N}&\cdots&\tfrac{1}{N}+a&\tfrac{2}{N}&\cdots&\tfrac{2}{N}&\cdots&\tfrac{N-1}{N}&\cdots&\tfrac{N-1}{N} \\
  \tfrac{1}{N}&\cdots&\tfrac{1}{N}&\tfrac{2}{N}+a&\cdots&\tfrac{2}{N}&\cdots&\tfrac{N-1}{N}&\cdots&\tfrac{N-1}{N} \\
  \cdots&\cdots&\cdots&\cdots&\cdots&\cdots&\cdots&\cdots&\cdots&\cdots \\
  \tfrac{1}{N}&\cdots&\tfrac{1}{N}&\tfrac{2}{N}&\cdots&\tfrac{2}{N}+a&\cdots&\tfrac{N-1}{N}&\cdots&\tfrac{N-1}{N} \\
 \cdots&\cdots&\cdots&\cdots&\cdots&\cdots&\cdots&\cdots&\cdots&\cdots \\
 \tfrac{1}{N}&\cdots&\tfrac{1}{N}&\tfrac{2}{N}&\cdots&\tfrac{2}{N}&\cdots&\tfrac{N-1}{N}+a&\cdots&\tfrac{N-1}{N}\\
  \cdots&\cdots&\cdots&\cdots&\cdots&\cdots&\cdots&\cdots&\cdots&\cdots \\
 \tfrac{1}{N}&\cdots&\tfrac{1}{N}&\tfrac{2}{N}&\cdots&\tfrac{2}{N}&\cdots&\tfrac{N-1}{N}&\cdots&\tfrac{N-1}{N}+a
\end{array}
\right].
$$
It is easy to check that $A$ is nonsingular for $a\neq 0$. Let
$\mathbf{A}$ be the matrix obtained from $A$ as follows: for an integer $j$ with $M_{i-1}+1\leq j\leq M_i$ for some $1\leq i\leq N-1$, the $j$-th row of $\mathbf{A}$ is obtained by multiplying the $j$-th row of $A$ by $$\frac{(\frac{i}{N})(2g-3+\sum_{i=1}^{N-1}n_i)!}{(a+\frac{i}{N})!\prod_{i=1}^{N-1}(\frac{i}{N})^{n_i}}.$$
Here $(a+\frac{i}{N})!:=\prod_{m=0}^a(m+\frac{i}{N}).$

The linear system
\begin{equation}\label{linearsystem}
\mathbf{A}\cdot \mathbf{c_{g}}=\mathbf{\Gamma_{g}}
\end{equation}
has a unique solution which represents $c_{j,g}$ as a linear combination of $\Gamma_{j,g}$'s for $1\leq j\leq \Sigma_{i}n_{i}$.

For integers $1\leq s\leq N-1$ and $1\leq j\leq\Sigma_{i}n_{i}$, we put $(k_{j}+\frac{s}{N})!=\prod_{m=0}^{k_j}(m+\frac{s}{N})$.

The vanishing of $z_{h,k_{1},\cdots,k_{\Sigma_{i}n_{i}}}^{k}(\mathbf{t})$ for $k\geq 1$ and $k_{1},\cdots,k_{\Sigma_{i}n_{i}}\geq 0$ yields the following theorem.

\begin{thm}\Label{curvevir}
We have
$$
\langle\widetilde{\tau}_{k_{1}}\cdots\widetilde{\tau}_{k_{\Sigma_{i}n_{i}}}\tau_{l_{1}}
\cdots\tau_{l_{n}}|\lambda_{r_{1}}\rangle^{h}_{g}=
\sum_{s=1}^{N-1}\sum_{j=M_{s-1}+1}^{M_{s}}\tfrac{(2g-3+n+\Sigma_{i}
n_{i})!(k_{j}+\frac{s}{N})}
{\prod_{j}l_{j}!\prod_{b=1}^{N-1}\prod_{j=M_{b-1}+1}^{M_{b}}(k_{j}+\frac{b}{N})!}c_{j,g}.
$$
\end{thm}

\begin{proof}
The vanishing of
$z_{h,k_{1},\cdots,k_{\Sigma_{i}n_{i}}}^{k}(\mathbf{t})$ gives the vanishing of its Taylor coefficients, $$\tfrac{1}{[1]_{0}^{k}}\partial_{t_{l_{1}}}\cdots \partial_{t_{l_{n}}}z_{h;k_{1},\cdots,k_{\Sigma_{i}n_{i}}}^{k}(0)=0.$$
An explicit calculation shows that $\tfrac{1}{[1]_{0}^{k}}\partial_{t_{l_{1}}}\cdots \partial_{t_{l_{n}}}z_{h;k_{1},\cdots,k_{\Sigma_{i}n_{i}}}^{k}(0)$ is the right side of the following
\begin{multline}\Label{recursion}
0=-\langle\tau_{k+1}\widetilde{\tau}_{k_{1}}\cdots\widetilde{\tau}_{k_{\Sigma_{i}n_{i}}}\tau_{l_{1}}
\cdots\tau_{l_{n}}|\lambda_{r_{1}}\rangle^{h}_{g}+\sum_{i=1}^{n}\tfrac{(l_{i}+k)!}{(l_{i}-1)!(k+1)!}\langle\widetilde{\tau}_{k_{1}}\cdots\widetilde{\tau}_{k_{\Sigma_{i}n_{i}}}\tau_{l_{1}}
\cdots\tau_{l_{i}+k}\cdots\tau_{l_{n}}|\lambda_{r_{1}}\rangle^{h}_{g}\\
+
\sum_{i=1}^{N-1}\sum_{j=M_{i-1}+1}^{M_{i}}\tfrac{(k_{j}+k+\frac{i}{N})!}{(k_{j}-1+\frac{i}{N})!(k+1)!}
\langle\widetilde{\tau}_{k_{1}}\cdots\widetilde{\tau}_{k_{M_{i-1}}}\widetilde{\tau}_{k_{M_{i-1}+1}}\cdots\widetilde{\tau}_{k_{j}+k}\cdots
\widetilde{\tau}_{k_{M_{i}}}\widetilde{\tau}_{k_{M_{i}+1}}\cdots\widetilde{\tau}_{k_{\Sigma_{i}n_{i}}}\tau_{l_{1}}
\cdots\tau_{l_{n}}|\lambda_{r_{1}}\rangle^{h}_{g}.
\end{multline}

We now solve the recursion (\ref{recursion}). The virtual dimension of $\overline{\mathcal{M}}_{g,n+1+\Sigma_{i}n_{i}}(\mathbb{P}(1,N),0,\mathbf{x})$ is
$$\text{vdim}=(1-g)(1-3)+n+\sum_{i=1}^{N-1}n_{i}+1-\sum_{i=1}^{N-1}n_{i}\frac{i}{N}.$$
If $\langle\tau_{k+1}\widetilde{\tau}_{k_{1}}\cdots\widetilde{\tau}_{k_{\Sigma_{i}n_{i}}}\tau_{l_{1}}\cdots\tau_{l_{n}}|\lambda_{r_{1}}\rangle^{h}_{g}\neq 0$, we have
$\text{vdim}=k+1+\sum_{i=1}^{n}l_{i}+\sum_{i=1}^{\Sigma_{i}n_{i}}k_{i}$. So

\begin{equation}\Label{virdim}
2g-2+n+\sum_{i=1}^{N-1}
n_{i}=\sum_{i=1}^{n}l_{i}+\sum_{i=1}^{N-1}\sum_{j=M_{i-1}+1}^{M_{i}}\Bigl(k_{j}+\tfrac{i}{N}\Bigr)+k.
\end{equation}

For an integer $r$, then there exists a unique  integer $s$ such that $M_{s-1}+1\leq r\leq M_{s}$, with $1\leq s\leq N-1$. 
Let $\mathbf{k}=(k_{1},\cdots,k_{\Sigma_{i}n_{i}})$ and 
$\mathbf{l}=(l_{1},\cdots,l_{n})$. Introduce 

\begin{equation}\label{special}
\Theta(\mathbf{k},\mathbf{l})_{r}:=\frac{(2g-3+n+\Sigma_{i}
n_{i})!(k_{r}+\frac{s}{N})}
{\prod_{j}l_{j}!\prod_{b=1}^{N-1}\prod_{j=M_{b-1}+1}^{M_{b}}(k_{j}+\frac{b}{N})!}.
\end{equation}

We claim that $\Theta(\mathbf{k},\mathbf{l})_{r}$ is a solution of the recursion (\ref{recursion}). To see this, write (\ref{virdim})
as
\begin{equation}\Label{virdimm}
2g-2+n+\sum_{i=1}^{N-1}
n_{i}=\sum_{i=1}^{n}l_{i}+\sum_{i=1}^{N-1}\sum_{j=M_{i-1}+1,j\neq r}^{M_{i}}\Bigl(k_{j}+\tfrac{i}{N}\Bigr)+(k_{r}+\frac{s}{N}+k).
\end{equation}
Multiply both sides of (\ref{virdimm}) by
$$
\frac{(2g-3+n+\Sigma_{i}n_{i})!(k_{r}+\tfrac{s}{N})}
{(k+1)!\prod_{j}l_{j}!\prod_{b=1}^{N-1}\prod_{j=M_{b-1}+1}^{M_{b}}(k_{j}+\frac{b}{N})!},
$$
we get
\begin{equation*}
\begin{split}
&\frac{(2g-2+n+\Sigma_{i}n_{i})!(k_{r}+\frac{s}{N})}
{(k+1)!\prod_{j}l_{j}!\prod_{b=1}^{N-1}\prod_{j=M_{b-1}+1}^{M_{b}}(k_{j}+\frac{b}{N})!}\\
&= \sum_{i=1}^{n}\frac{(l_{i}+k)!}{(l_{i}-1)!(k+1)!}\frac{(2g-3+n+\Sigma_{i}n_{i})!(k_{r}+\frac{s}{N})}
{l_{1}!\cdots (l_{i}+k)!\cdots l_{n}!\prod_{b=1}^{N-1}\prod_{j=M_{b-1}+1}^{M_{b}}(k_{j}+\frac{b}{N})!} \\
&+\sum_{i=1}^{N-1}\frac{(k_{j}+k+\frac{i}{N})!}{(k_{j}-1+\frac{i}{N})!(k+1)!}\cdot
\sum_{j=M_{i-1}+1,j\neq r}^{M_{i}} \tfrac{(2g-3+n+\Sigma_{i}n_{i})!(k_{r}+\frac{s}{N})}
{\prod_{j}l_{j}!\prod_{b=1,b\neq i}^{N-1}\prod_{j=M_{b-1}+1}^{M_{b}}(k_{j}+\frac{b}{N})!(k_{M_{i-1}+1}+\frac{i}{N})!\cdots (k_{j}+k+\frac{i}{N})!\cdots(k_{M_{i}}+\frac{i}{N})!}\\
&+\frac{(k_{r}+k+\frac{s}{N})!}{(k_{r}-1+\frac{s}{N})!(k+1)!}\cdot \tfrac{(2g-3+n+\Sigma_{i}n_{i})!(k_{r}+k+\frac{s}{N})}{\prod_{j}l_{j}!\prod_{b=1,b\neq s}^{N-1}\prod_{j=M_{b-1}+1}^{M_{b}}(k_{j}+\frac{b}{N})!(k_{M_{s-1}+1}+\frac{s}{N})!\cdots (k_{r}+k+\frac{s}{N})!\cdots(k_{M_{s}}+\frac{i}{N})!}.
\end{split}
\end{equation*}
It is straightforward to see that this is the recursion (\ref{recursion}). 

Suppose that $\langle\widetilde{\tau}_{k_{1}}\cdots\widetilde{\tau}_{k_{\Sigma_{i}n_{i}}}\tau_{l_{1}}
\cdots\tau_{l_{n}}|\lambda_{r_{1}}\rangle^{h}_{g}$ is of the form $\sum_r c_{r,g}\Theta(\mathbf{k},\mathbf{l})_r$. Then by considering special values of $\mathbf{k}, \mathbf{l}$, we find that the coefficients $c_{r,g}$ are uniquely determined by the linear system (\ref{linearsystem}). The result follows.
\end{proof}

The initial values (\ref{constant1}) in Theorem \ref{curvevir} and more generally the following Hurwitz-Hodge integrals 
$$\int_{\overline{\M}_{g,\sum n_{i}}(\mathcal{B}\mathbb{Z}_{N})}\lambda_{r_{1}-i}\psi^{a+i}, \quad i\geq 0,$$
are computed in \cite{JPT}. Details can be found there.

\section{Virasoro constraints for weighted projective planes}\label{weighted_P2}

Let $\X=\mathbb{P}(1,1,N)$. The Chen-Ruan orbifold cohomology $H^{*}_{CR}(\mathbb{P}(1,1,N))$
has the following natural generators:
\begin{equation*}
\begin{split}
& 1\in H^{0}_{CR}(\mathbb{P}(1,1,N)), \,\xi\in H^{2}_{CR}(\mathbb{P}(1,1,N)), \, [\mathcal{X}]\in H^{4}_{CR}(\mathbb{P}(1,1,N)),\\
& \gamma_{j}\in H^{\tfrac{4j}{N}}_{CR}(\mathbb{P}(1,1,N)) \, \text{ for }1\leq j\leq N-1.
\end{split}
\end{equation*}

The following coordinates for the corresponding descendants will be used:
\begin{center}
\begin{tabular}{lp{0.75\textwidth}}
coordinate & \quad descendant\\
$t_i, \, i\geq 0$ & \quad $\tau_i(1)$\\
$s_{i}, \, i\geq 0$ & \quad $\tau_{i}(\xi)$\\
$r_{i}, \, i\geq 0$ & \quad $\tau_{i}([\mathcal{X}])$\\
$\alpha^{j}_{i}, \, i\geq 0$ & \quad $\tau_{i}(\gamma_{j})$ for $1\leq j\leq N-1$.
\end{tabular}
\end{center}

Let $\mathcal{F}_{\mathcal{X}}^{g,0}$ be the degree zero, genus $g$
orbifold Gromov-Witten potential of $\mathcal{X}$. Let
$$\mathcal{D}_{\mathcal{X}}^{0}=\exp\left(\sum_{g\geq 0}\hbar^{g-1}\mathcal{F}_{\mathcal{X}}^{g,0}\right).$$ Let $c$ be the number such that $c_{1}(\mathcal{X})=c\cdot\xi$.
\begin{thm}
We have
\begin{align*}
\frac{\mathcal{L}_{k}\mathcal{D}_{\mathcal{X}}^{0}}{\mathcal{D}_{\mathcal{X}}^{0}}=
\sum_{g=0}^{\infty}\hbar^{g-1}(-1)^{g}\left(|c|^{2}x_{g}^{k}(\mathbf{t})-\sum_{l=0}^{\infty}c\cdot s_{l}y_{g,l}^{k}(\mathbf{t})\right)+ \frac{1}{\hbar}w(r,s,t) 
+\sum_{g=0}^{\infty}\hbar^{g-1}(-1)^{g}\sum_{h=0}^{\infty} u_{g,l}^{h,k}(\mathbf{t}) \\
+\sum_{g=0}^{\infty}\hbar^{g-1}(-1)^{g}\sum_{h=0}^{\infty}~\sum_{k_{1},\cdots,k_{\Sigma_{i}n_{i}}> 0}\alpha^{1}_{k_{1}}\cdots\alpha^{N-1}_{k_{\Sigma_{i}
n_{i}}}\cdot
z_{h;k_{1},\cdots,k_{\Sigma_{i}n_{i}}}^{k}(\mathbf{t}),
\end{align*}
where 
\begin{equation*}
\begin{split}
x_{g}^{k}(\mathbf{t})=
&-[\tfrac{1}{2}]_{0}^{k}\langle\langle\tau_{k+1}|\lambda_{g}\lambda_{g-2}\rangle\rangle_{g}
-\sum_{m=0}^{\infty}[m-\tfrac{1}{2}]_{0}^{k}t_{m}\langle\langle\tau_{k+m}|\lambda_{g}\lambda_{g-2}\rangle\rangle_{g}\\
&+\sum_{m=0}^{k-1}(-1)^{m+1}\Bigl( [-m-\tfrac{3}{2}]_{0}^{k}\langle\langle\tau_{m}\rangle\rangle \langle\langle\tau_{k-m-1}|\lambda_{g}\lambda_{g-2}\rangle\rangle_{g} \\
&+\tfrac{1}{2}[-m-\tfrac{1}{2}]_{0}^{k}\sum_{g=g_{1}+g_{2}}\langle\langle\tau_{m}|\lambda_{g_{1}}\lambda_{g_{1}-2}\rangle\rangle_{g_{1}} \langle\langle\tau_{k-m-1}|\lambda_{g_{2}}\lambda_{g_{2}-2}\rangle\rangle_{g_{2}}\Bigr) \\
&+[\tfrac{1}{2}]_{1}^{k}\langle\langle \tau_{k}|\lambda_{g}\lambda_{g-1}\rangle\rangle_{g}+ \sum_{m=0}^{\infty}t_{m}[m-\tfrac{1}{2}]_{1}^{k}\langle\langle\tau_{k+m-1}|\lambda_{g}\lambda_{g-1}\rangle\rangle_{g} \\
&+ \sum_{m=0}^{k-2}(-1)^{m+1}[-m-\tfrac{3}{2}]_{1}^{k}\langle\langle\tau_{m}\rangle\rangle_{0}\langle\langle\tau_{k-m-2}|\lambda_{g}\lambda_{g-1}\rangle\rangle_{g},
\end{split}
\end{equation*}

\begin{equation*}
\begin{split}
y_{g,l}^{k}(\mathbf{t})=
-[\tfrac{1}{2}]_{0}^{k}\langle\langle\tau_{k+1}\tau_{l}|\lambda_{g}\lambda_{g-1}\rangle\rangle_{g}
+\sum_{m=0}^{\infty}[m-\tfrac{1}{2}]_{0}^{k}t_{m}\langle\langle\tau_{k+m}\tau_{l}|\lambda_{g}\lambda_{g-1}\rangle\rangle_{g}
+[l+\tfrac{1}{2}]_{0}^{k}\langle\langle\tau_{k+l}|\lambda_{g}\lambda_{g-1}\rangle\rangle_{g}\\
+ \sum_{m=0}^{k-1}(-1)^{m+1}\Bigl([-m-\tfrac{3}{2}]_{0}^{k}\langle\langle\tau_{m}\rangle\rangle_{0}\langle\langle\tau_{k-m-1}\tau_{l}|\lambda_{g}\lambda_{g-1}\rangle\rangle_{g}\\ 
+[-m-\tfrac{1}{2}]_{0}^{k}\langle\langle\tau_{m}\tau_{l}\rangle\rangle_{0}\langle\langle\tau_{k-m-1}|\lambda_{g}\lambda_{g-1}\rangle\rangle_{g}\Bigr),
\end{split}
\end{equation*}
$w(r,s,t)$ coincides with a similar term in \cite{GP}, and
$$u_{g,l}^{h,k}(\mathbf{t})=-[\tfrac{1}{2}]_{0}^{k}\langle\langle\tau_{k+1}|\lambda_{r_{1}}^{2}\rangle\rangle_{g}^{h}+
\sum_{m=1}^{\infty}[m-\tfrac{1}{2}]_{0}^{k}t_{m}\langle\langle\tau_{k+m}|\lambda_{r_{1}}^{2}\rangle\rangle_{g}^{h},$$

\begin{multline*}
z_{h;k_{1},\cdots,k_{\Sigma_{i}n_{i}}}^{k}(\mathbf{t})=
-[\tfrac{1}{2}]_{0}^{k}\langle\langle\tau_{k+1}\widetilde{\tau}_{k_{1}}\cdots
\widetilde{\tau}_{k_{\Sigma_{i}n_{i}}}|\lambda^{2}_{r_{1}}\rangle\rangle^{h}_{g}
+\sum_{m=1}^{\infty}[m-\tfrac{1}{2}]_{0}^{k}t_{m}\langle\langle\tau_{k+m}\widetilde{\tau}_{k_{1}}\cdots \widetilde{\tau}_{k_{\Sigma_{i}n_{i}}}|\lambda^{2}_{r_{1}}\rangle\rangle^{h}_{g}\\
+\sum_{i=1}^{N-1}\sum_{j=M_{i-1}+1}^{M_{i}}[k_{j}+\tfrac{2i}{N}-\tfrac{1}{2}]_{0}^{k}
\langle\langle\widetilde{\tau}_{k_{1}}\cdots\widetilde{\tau}_{k_{n_{i-1}}}\widetilde{\tau}_{k_{n_{i-1}+1}}\cdots\widetilde{\tau}_{k+k_{j}}\cdots\widetilde{\tau}_{k_{n_{i}}}\widetilde{\tau}_{k_{n_{i}+1}}
\cdots\widetilde{\tau}_{k_{\Sigma_{i}n_{i}}}|\lambda^{2}_{r_{1}}\rangle\rangle^{h}_{g}.
\end{multline*}
\end{thm}

\begin{proof}
The Virasoro operator $\mathcal{L}_{k}$, $k>0$ is given by
\begin{align*}
\mathcal{L}_{k}&= -[\tfrac{1}{2}]_{0}^{k}\partial_{t_{k+1}}+\sum_{m=0}^{\infty}\Bigl([m-\tfrac{1}{2}]_{0}^{k}t_{m}\partial_{t_{k+m}}
+[m+\tfrac{1}{2}]_{0}^{k}s_{m}\partial_{s_{k+m}}+[m+\tfrac{3}{2}]_{0}^{k}r_{m}\partial_{r_{k+m}}\\
&+\sum_{i=1}^{N-1}[m+\tfrac{2i}{N}-\tfrac{1}{2}]_{0}^{k}\alpha_{m}^{i}\partial_{\alpha_{k+m}^{i}}\Bigr)
+\hbar\sum_{m=0}^{k-1}(-1)^{m+1}\Bigl([-m-\tfrac{3}{2}]_{1}^{k}\partial_{r_{m}}\partial_{t_{k-m-1}}+\\
&\tfrac{1}{2}[-m-\tfrac{1}{2}]_{1}^{k}\partial_{s_{m}}\partial_{s_{k-m-1}}
+ \sum_{i=1}^{N-1}[-m-\tfrac{2i}{N}+\tfrac{1}{2}]_{1}^{k}\partial_{\alpha^{i}_{m}}\partial_{\alpha^{N-i}_{k-m-1}}\Bigr) + \\
&\mathbf{c}\cdot\Bigl(-[\tfrac{1}{2}]_{1}^{k}\partial_{s_{k}}+\sum_{m=0}^{\infty}[m-\tfrac{1}{2}]_{1}^{k}t_{m}\partial_{s_{k+m-1}} 
+[m+\tfrac{1}{2}]_{1}^{k}s_{m}\partial_{r_{k+m-1}}\\
&+\hbar\sum_{m=0}^{k-2}(-1)^{m+1}[-m-\tfrac{3}{2}]_{1}^{k}\partial_{r_{m}}\partial_{s_{k-m-2}}
\Bigr)
+|\mathbf{c}|^{2}\cdot\Bigl(-[\tfrac{1}{2}]_{2}^{k}\partial_{r_{k-1}}
 +\sum_{m=0}^{\infty}[m-\tfrac{1}{2}]_{2}^{k}t_{m}\partial_{r_{k+m-2}}\\
&+\frac{\hbar}{2}\sum_{m=0}^{k-3}(-1)^{m+1}[-m-\tfrac{3}{2}]_{2}^{k}\partial_{r_{m}}\partial_{r_{k-m-3}} \Bigr)+\frac{\delta_{k1}}{2\hbar}t_{0}^{2},
\end{align*}
where $\mathbf{c}$ satisfies $c_{1}(\mathcal{X})=\mathbf{c}\xi$.

Similar to the one-dimensional case, contributions from components which parametrize maps with trivial monodromies need extra care. Such a component in the moduli space of degree $0$ genus $g$ $n$-pointed stable maps is a product $\overline{\M}_{g,n}\times \mathbb{P}(1,1,N)$. As in \cite{GP} the obstruction bundle over this component is identified with $T_{\mathbb{P}(1,1,N)}\boxtimes \mathbb{E}^{\vee}$ and contributions from this component are integrals against the Euler class 
$$e(T_{\mathbb{P}(1,1,N)}\boxtimes \mathbb{E}^{\vee})=-c_{1}(\mathbb{P}(1,1,N))\lambda_{g}\lambda_{g-1}-c_{1}(\mathbb{P}(1,1,N))^{2}\lambda_{g}\lambda_{g-2}.$$

Thus the degree zero orbifold Gromov-Witten potential is:
\begin{multline*}
\mathcal{D}_{\mathcal{X}}^{0}=
\exp\Bigl(\frac{1}{\hbar}\sum_{m=0}^{\infty}r_{m}\langle\langle\tau_{m}\rangle\rangle_{0}
+\frac{1}{\hbar}\sum_{l,m}\frac{1}{2}s_{l}s_{m}\langle\langle \tau_{l}\tau_{m}\rangle\rangle_{0}
+|c|^{2}\sum_{g=1}^{\infty}\hbar^{g-1}\langle\langle |\lambda_{g}\lambda_{g-2}\rangle\rangle_{g}\\-
\sum_{g=1}^{\infty}\hbar^{g-1}\sum_{m=0}^{\infty}c\cdot s_{m}\langle\langle\tau_{m}|\lambda_{g}\lambda_{g-1}\rangle\rangle_{g}
+\sum_{g=0}^{\infty}(-1)^{g}\hbar^{g-1}\sum_{h=0}^{\infty}~\langle\langle~|\lambda_{r_{1}}^{2}\rangle\rangle_{g}^{h}\\
+\sum_{g=0}^{\infty}(-1)^{g}\hbar^{g-1}\sum_{h=0}^{\infty}~
\sum_{k_{1},\cdots,k_{\Sigma_{i}n_{i}}> 0}\alpha^{1}_{k_{1}}
\cdots\alpha^{N-1}_{k_{\Sigma_{i}
n_{i}}}\langle\langle\widetilde{\tau}_{k_{1}}\cdots
\widetilde{\tau}_{k_{\Sigma_{i}n_{i}}}|\lambda^{2}_{r_{1}}\rangle\rangle^{h}_{g}\Bigr).
\end{multline*}
So applying the operator $\mathcal{L}_{k}$ to
$\mathcal{D}_{\mathcal{X}}^{0}$ we obtain the result.
\end{proof}

The degree zero Virasoro constraints for $\mathbb{P}(1,1,N)$ is equivalent to the vanishing of $x_{g}^{k}(\mathbf{t})$, $y_{g,l}^{k}(\mathbf{t})$, $w(r,s,t)$ and $u_{g,l}^{h,k}(\mathbf{t})$ and  $z_{h,k_{1},\cdots,k_{\Sigma_{i}n_{i}}}^{k}(\mathbf{t})$. As explained in \cite{GP}, the formula $w(r,s,t)$ is very complicated. To get formula of Hurwitz-Hodge integrals, it is not necessary to write $w(r,s,t)$ down.

\subsubsection{\textbf{Vanishing of $u_{g,l}^{h,k}(\mathbf{t})$}.}

\begin{thm}\Label{surfacevir-1}
We have
$$
\langle\tau_{l_{1}}
\cdots\tau_{l_{n}}|\lambda^{2}_{r_{1}}\rangle^{h}_{g}=
\frac{(2g+n-3)!(2g-1)!!}{(2g-1)! (2l_{1}-1)!!\cdots (2l_{n}-1)!!}\langle\tau_{g}|\lambda^{2}_{r_{1}}\rangle_g^h. 
$$
\end{thm}

\begin{proof}
The vanishing of
$u_{g,l}^{h,k}(\mathbf{t})$ gives the vanishing of its Taylor coefficients, $$\tfrac{1}{[\tfrac{1}{2}]_{0}^{k}}\partial_{t_{l_{1}}}\cdots \partial_{t_{l_{n}}}u_{g,l}^{h,k}(0)=0.$$
An explicit calculation shows that $\tfrac{1}{[\tfrac{1}{2}]_{0}^{k}}\partial_{t_{l_{1}}}\cdots \partial_{t_{l_{n}}}u_{g,l}^{h,k}(0)$ is the right side of the following
\begin{equation}\Label{recursion-sur-1}
0=-\langle\tau_{k+1}\tau_{l_{1}}
\cdots\tau_{l_{n}}|\lambda^{2}_{r_{1}}\rangle^{h}_{g}+\sum_{i=1}^{n}\tfrac{(2l_{i}+2k-1)!!}{(2k+1)!!(2l_{i}-1)!!}\langle\tau_{l_{1}}
\cdots\tau_{l_{i}+k}\cdots\tau_{l_{n}}|\lambda^{2}_{r_{1}}\rangle^{h}_{g}
\end{equation}
The result follows by solving the recursion (\ref{recursion-sur-1}). 
\end{proof}

Exact evaluations of initial values in Theorem \ref{surfacevir-1} seem to be unavailable at the moment. 

\subsubsection{\textbf{Vanishing of $z_{h,k_{1},\cdots,k_{\Sigma_{i}n_{i}}}^{k}(\mathbf{t})$}.}
Let
\begin{equation}\Label{constant2}
\Gamma_{j,g}=\langle\widetilde{\tau}_{a}|\lambda^{2}_{r_{1}}\rangle^{h}_{g}=\int_{\Mbar_{g,\sum n_{i}}(\mathcal{B}\mathbb{Z}_{N})}\psi^{a}\lambda^{2}_{r_{1}},\quad \text{where } a:=g-1+\sum_{i=1}^{N-1}n_{i}-\sum_{i=1}^{N-1}\frac{2i}{N}n_{i}.
\end{equation}

Again write $\mathbf{\Gamma_{g}}:=(\Gamma_{j,g})_{1\leq j\leq \Sigma_{i}n_{i}}$ as a column vector and let $\mathbf{c_{g}}:=(c_{j,g})_{1\leq j\leq \Sigma_{i}n_{i}}$ be another column vector. Let the index $i$ vary from $1$ to $N-1$ and define a $\Sigma_{i}n_{i}\times\Sigma_{i}n_{i}$ square matrix $A=(a_{st})$ by
\begin{equation}\Label{matrix2}
a_{st}:=\begin{cases}\frac{2i}{N}+a&\text{if
}~ M_{i-1}<s=t\leq M_{i}\,;\\
\frac{2i}{N}&\text{if}~M_{i-1}<t\leq M_{i}~\text{and}~s\neq t\,.\end{cases}
\end{equation}

The matrix can be written as follows:
$$
A=\left[
\begin{array}{cccccccccc}
  \tfrac{2}{N}+a&\cdots&\tfrac{2}{N}&\tfrac{4}{N}&\cdots&\tfrac{4}{N}&\cdots&\tfrac{2(N-1)}{N}&\cdots&\tfrac{2(N-1)}{N} \\
  \tfrac{2}{N}&\cdots&\tfrac{2}{N}&\tfrac{4}{N}&\cdots&\tfrac{4}{N}&\cdots&\tfrac{2(N-1)}{N}&\cdots&\tfrac{2(N-1)}{N} \\
  \cdots&\cdots&\cdots&\cdots&\cdots&\cdots&\cdots&\cdots&\cdots&\cdots\\
  \tfrac{2}{N}&\cdots&\tfrac{2}{N}+a&\tfrac{4}{N}&\cdots&\tfrac{4}{N}&\cdots&\tfrac{2(N-1)}{N}&\cdots&\tfrac{2(N-1)}{N} \\
  \tfrac{2}{N}&\cdots&\tfrac{2}{N}&\tfrac{4}{N}+a&\cdots&\tfrac{4}{N}&\cdots&\tfrac{2(N-1)}{N}&\cdots&\tfrac{2(N-1)}{N} \\
  \cdots&\cdots&\cdots&\cdots&\cdots&\cdots&\cdots&\cdots&\cdots&\cdots \\
  \tfrac{2}{N}&\cdots&\tfrac{2}{N}&\tfrac{4}{N}&\cdots&\tfrac{4}{N}+a&\cdots&\tfrac{2(N-1)}{N}&\cdots&\tfrac{2(N-1)}{N} \\
 \cdots&\cdots&\cdots&\cdots&\cdots&\cdots&\cdots&\cdots&\cdots&\cdots \\
 \tfrac{2}{N}&\cdots&\tfrac{2}{N}&\tfrac{4}{N}&\cdots&\tfrac{4}{N}&\cdots&\tfrac{2(N-1)}{N}+a&\cdots&\tfrac{2(N-1)}{N}\\
  \cdots&\cdots&\cdots&\cdots&\cdots&\cdots&\cdots&\cdots&\cdots&\cdots \\
 \tfrac{2}{N}&\cdots&\tfrac{2}{N}&\tfrac{4}{N}&\cdots&\tfrac{4}{N}&\cdots&\tfrac{2(N-1)}{N}&\cdots&\tfrac{2(N-1)}{N}+a
\end{array}
\right].
$$
It is easy to check that $A$ is nonsingular for $a\neq 0$. Let $\mathbf{A}$ be the matrix obtained from $A$ as follows: for an integer $j$ with $M_{i-1}+1\leq j\leq M_i$ for some $1\leq i\leq N-1$, the $j$-th row of $\mathbf{A}$ is obtained by multiplying the $j$-th row of $A$ by
$$\frac{(g+\frac{1}{2}(\sum_{i}n_i-3))!(\frac{2i}{N}-\frac{1}{2})}{(a+\frac{2i}{N}-\frac{1}{2})!\prod_{i=1}^{N-1}(\frac{2i}{N}-\frac{1}{2})^{n_i}}.$$

The linear system
\begin{equation}\label{linearsystem2}
\mathbf{A}\cdot \mathbf{c_{g}}=\mathbf{\Gamma_{g}}
\end{equation}
has a unique solution which represents $c_{j,g}$ as a linear combination of $\Gamma_{j,g}$'s for $1\leq j\leq \Sigma_{i}n_{i}$. 

For integers $1\leq s\leq N-1$ and $1\leq j\leq\Sigma_{i}n_{i}$, let $(k_{j}+\frac{2s}{N})!=\frac{2s}{N}\cdot(1+\frac{2s}{N})\cdots(k_{j}+\frac{2s}{N})$.

The vanishing of $z_{h,k_{1},\cdots,k_{\Sigma_{i}n_{i}}}^{k}(\mathbf{t})$ for $k\geq 1$ and $k_{1},\cdots,k_{\Sigma_{i}n_{i}}\geq 0$ yields the following result. 

\begin{thm}\Label{surfacevir}
We have
$$
\langle\widetilde{\tau}_{k_{1}}\cdots\widetilde{\tau}_{k_{\Sigma_{i}n_{i}}}\tau_{l_{1}}
\cdots\tau_{l_{n}}|\lambda^{2}_{r_{1}}\rangle^{h}_{g}=
\sum_{s=1}^{N-1}\sum_{j=M_{s-1}+1}^{M_{s}}\tfrac{(g+\frac{1}{2}(n+\Sigma_{i}
n_{i}-3))!(k_{j}-\frac{1}{2}+\frac{2s}{N})}
{\prod_{j}(l_{j}-\frac{1}{2})!\prod_{b=1}^{N-1}\prod_{j=M_{b-1}+1}^{M_{b}}(k_{j}-\frac{1}{2}+\frac{2b}{N})!}c_{j,g}.
$$
\end{thm}

\begin{proof}
Again we consider the following recursion given by  $\tfrac{1}{[\frac{1}{2}]_{0}^{k}}\partial_{t_{l_{1}}}\cdots \partial_{t_{l_{n}}}z_{h;k_{1},\cdots,k_{\Sigma_{i}n_{i}}}^{k}(0)=0$:

\begin{multline}\Label{recursion2}
0=-\langle\tau_{k+1}\widetilde{\tau}_{k_{1}}\cdots\widetilde{\tau}_{k_{\Sigma_{i}n_{i}}}\tau_{l_{1}}
\cdots\tau_{l_{n}}|\lambda^{2}_{r_{1}}\rangle^{h}_{g}
+\sum_{i=1}^{n}\tfrac{[l_{i}-\frac{1}{2}]^{k}_{0}}{[\frac{1}{2}]_{0}^{k}}\langle\widetilde{\tau}_{k_{1}}\cdots\widetilde{\tau}_{k_{\Sigma_{i}n_{i}}}\tau_{l_{1}}
\cdots\tau_{l_{i}+k}\cdots\tau_{l_{n}}|\lambda^{2}_{r_{1}}\rangle^{h}_{g}\\
+\sum_{i=1}^{N-1}\sum_{j=M_{i-1}+1}^{M_{i}}\tfrac{[k_{j}+\frac{2i}{N}-\frac{1}{2}]_{0}^{k}}{[\frac{1}{2}]_{0}^{k}}
\langle\widetilde{\tau}_{k_{1}}\cdots\widetilde{\tau}_{k_{n_{i-1}}}\widetilde{\tau}_{k_{n_{i-1}+1}}\cdots\widetilde{\tau}_{k_{j}+k}\cdots
\widetilde{\tau}_{k_{n_{i}}}\widetilde{\tau}_{k_{n_{i}+1}}\cdots\widetilde{\tau}_{k_{\Sigma_{i}n_{i}}}\tau_{l_{1}}
\cdots\tau_{l_{n}}|\lambda^{2}_{r_{1}}\rangle^{h}_{g}.
\end{multline}
Virtual dimension constraints for orbifold Gromov-Witten invariants of $\mathbb{P}(1,1,N)$ gives 
\begin{equation}\label{virdim2}
g+\frac{1}{2}\Bigl(n+\sum_{i=1}^{N-1}
n_{i}-2\Bigr)
=\sum_{i=1}^{n}\Bigl(l_{i}-\tfrac{1}{2}\Bigr)+\sum_{i=1}^{N-1}\sum_{j=M_{i-1}+1}^{M_{i}}\Bigl(k_{j}-\tfrac{1}{2}+\tfrac{2i}{N}\Bigr)+k.
\end{equation}
So from (\ref{linearsystem2}), (\ref{recursion2}) and
(\ref{virdim2}), using the same method as in the proof of Theorem
\ref{curvevir} we finish the proof.
\end{proof}

Again, initial values in Theorem \ref{surfacevir} remain to be computed explicitly.

\section{Degree zero Virasoro constraints for threefolds}\label{threefolds}

The Virasoro constraints in degree zero for threefolds do not give anything new:
all descendent invariants are reduced to primary ones by string and dilaton equations. This follows from a dimension argument already discussed in \cite{JT}.

Degree zero primary invariants of $\mathbb{P}(1,1,1,N)$ are closely related to Gromov-Witten invariants of $[\mathbb{C}^3/\zz_N]$, the stack defined by the weight $(1,1,1)$ action of $\zz_N$ on $\mathbb{C}^3$. The case $N=3$ is of special interest since this action is Calabi-Yau. In this case primary invariants with stacky insertions are written as Hurwitz-Hodge integrals 
\begin{equation*}
\int_{\overline{\mathcal{M}}_{g,n_1+n_2}(\mathcal{B}\mathbb{Z}_{3})}\lambda^{3}_{r_{1}},
\end{equation*}
where $n_{1},n_{2}$ represent $n_{1}$ stacky points of type $\omega$
and $n_{2}$ stacky points of type $\overline{\omega}$.

Some of these integrals have been predicted in physics \cite{ABK}. For example,
$$\int_{\overline{\mathcal{M}}_{1,
\omega,\omega,\omega}(\mathcal{B}\mathbb{Z}_{3})}\lambda^{3}_{r_{1}}=0,\quad \int_{\overline{\mathcal{M}}_{1,\omega^{6}}(\mathcal{B}\mathbb{Z}_{3})}\lambda^{3}_{r_{1}}=\frac{1}{3^{5}}.$$

In genus zero, the generating function of these integrals is computed in \cite{ccit}, \cite{ccit2}, \cite{BaCa}, \cite{CaC}. Mathematical calculations of examples in higher genus are found in \cite{BoC}.


\end{document}